\documentclass[12pt,notitlepage]{amsart}
\usepackage{amsmath}     
\usepackage{amscd}

\advance\voffset by -0.4 truein
\advance\hoffset by -0.4 truein
\newcommand{\m}{{\mathfrak m}}

 \font\tenbi=cmmi14
 \font\sevenbi=cmmi10 \font\fivebi=cmmi7
 \newfam\bifam \def\bi{\fam\bifam} \textfont\bifam=\tenbi
 \scriptfont\bifam=\sevenbi \scriptscriptfont\bifam=\fivebi
 \mathchardef\variablemega="7121 \def\bigomega{{\bi\variablemega}}

\newcommand{\Ocal}{\mathcal{O}}

\newcommand{\wt}{\widetilde}
\newcommand{\ol}{\overline}
\renewcommand{\:}{\colon}
\newcommand{\IP}{\text{\bf P}_{\hskip-0.1cm k}}

\newcommand{\w}{\bigomega}
\renewcommand{\deg}{\mathrm{deg}\,}
\renewcommand{\d}{\partial}

\newtheorem{theorem}{Theorem}
\newtheorem{corollary}[theorem]{Corollary}
\newtheorem{proposition}[theorem]{Proposition}
\theoremstyle{definition}
\newtheorem{remark}[theorem]{Remark}
\newtheorem{subsct}[theorem]{}
\newtheorem{example}[theorem]{Example}
\theoremstyle{plain}
 
\font\smallrm=cmr8
\font\smallsc=cmcsc10
\font\smallsl=cmsl10

\begin{document}
\author
[{\smallrm EDUARDO ESTEVES}]
{Eduardo Esteves}
\title
[{\smallrm THE REGULARITY OF AN INTEGRAL VARIETY}]
{The Castelnuovo-Mumford regularity of \\ 
an integral variety of a vector field \\
on projective space}
\thanks{Supported by PRONEX, Conv\^enio 41/96/0883/00, CNPq, 
Proc. 300004/95-8, and FAPERJ, Proc. E-26/170.418/2000-APQ1}
\begin{abstract}
The Castelnuovo-Mumford regularity $r$ of a 
variety $V\subseteq\text{\bf P}^n_{\hskip-0.1cm\text{\bf C}}$ 
is an upper bound for the degrees of the hypersurfaces necessary to cut out 
$V$. In this note we give a bound for $r$ when $V$ is left invariant 
by a vector field on $\text{\bf P}^n_{\hskip-0.1cm\text{\bf C}}$. 
More precisely, assume $V$ is arithmetically Cohen-Macaulay, for instance, 
a complete intersection. Assume as well that $V$ projects to a 
normal-crossings hypersurface, which is the case when $V$ is a curve with 
at most ordinary nodes. Then 
we show that $r\leq m+s+1$, where $s$ is the dimension of $V$ and $m$ is 
the degree of the vector field. Our method consists of using first central 
projections to reduce the problem to when $V$ is a hypersurface, and then 
bounds given by Brunella and Mendes. 
\end{abstract}
\maketitle

\section{Introduction}
\label{section1}
Let $k$ be an algebraically closed field of characteristic $p\geq 0$, and 
$\IP^n$ the $n$-dimensional projective space over $k$. Let $X$ be 
an algebraic vector field on $\IP^n$. 
A subscheme $V\subseteq\IP^n$ is said to be left invariant by $X$, or 
an integral subscheme of $X$, 
if $X$ restricts to a vector 
field on $V$. 
For $k=\text{\bf C}$, the vector field $X$ induces a holomorphic 
flow away from its singular points. One might ask whether there are 
complex, compact subvarieties of 
$\text{\bf P}_{\hskip-0.1cm\text{\bf C}}^n$ where the 
flow stays confined. These subvarieties are algebraic by 
Chow's theorem, and are left invariant by $X$.

For $k=\text{\bf C}$ and $n=2$, H.~Poincar\'e \cite{P} considered the 
following question: Can we give a bound for the degrees of plane curves 
left invariant by $X$? If a bound were given, then we could 
restrict the search for such curves over a finite-dimensional vector space, 
that of polynomials of degree at most the given bound. In addition, it 
would be interesting if the bound depended only on the unique numerical 
invariant of $X$, its degree. 
For general $k$ and $n$, the degree of $X$ is the number of points on a 
general hyperplane $H\subseteq\IP^n$ at which 
the direction given by $X$ lies in $H$. Let $m$ denote the degree of $X$.

It's impossible to find the bound asked by Poincar\'e for $k=\text{\bf C}$ 
and $n=2$ only in terms of $m$; see Remark~9. 
However, given a reduced curve 
$C\subseteq\text{\bf P}^2_{\hskip-0.1cm\text{\bf C}}$ 
of degree $d$ left invariant by $X$, bounds for $d$ in terms of 
$m$ were found under certain conditions on $X$ or $C$. 
For instance, if the singularities of $X$ are non-diacritical, 
then M.~Carnicer showed in 
\cite{C} that $d\leq m+2$. 
If $C$ has at most ordinary nodes for singularities, D. Cerveau and 
A. Lins Neto showed in \cite{CL} that $d\leq m+2$ as well, with equality 
only if $C$ is reducible. 
If $C$ has worse singularities or $X$ is diacritical, the above inequality
does not hold. For these cases, weaker inequalities were given by 
A.~Campillo and Carnicer in \cite{CC}. 

For $k=\text{\bf C}$ but general $n$, 
M. Soares showed that the inequality $d\leq m+1$ holds for the 
degree $d$ of a smooth hypersurface left invariant by $X$; see \cite{So1}. 
Recently, M. Brunella and L.G. Mendes showed that a hypersurface 
of $\text{\bf P}^n_{\hskip-0.1cm\text{\bf C}}$ 
left invariant by $X$ has degree at most $m+n$ if it has at most 
normal-crossings singularities; see \cite{BM}. 
(More generally, they bounded the degrees of solutions to Pfaff equations 
on complex, projective manifolds $M$ with $\text{Pic}(M)\cong\mathbb Z$.) 

Further results are obtained for complete intersections. 
If $C\subseteq\text{\bf P}^n_{\hskip-0.1cm\text{\bf C}}$ 
is a complete intersection of 
hypersurfaces of degrees $d_1,\dots,d_{n-1}$, and is left invariant by $X$, 
a quite simple inequality was obtained by Soares in \cite{So2}:
\[
d_1+\cdots+d_{n-1}\leq m+n-1\quad\text{if $C$ is smooth.}
\]
If $C$ has at most ordinary nodes for singularities, Campillo, Carnicer 
and J. Garc\'\i a de la Fuente showed in \cite{CCF} that 
$\sum d_i\leq m+n$. They gave also inequalities for when $C$ has worse 
singularities.

In the present note we give a direct algebraic 
proof of Soares' result, which is essentially an observation 
by O. Zariski; see Theorem 8. Furthermore, the proof yields a 
characterization of all vector fields that leave invariant a given smooth 
hypersurface; see Theorem~7. 

Most importantly, we treat the case $n\geq 3$. 
Following up on Poincar\'e's original 
intention, we want to restrict to a finite-dimensional 
space the search for a subscheme $V\subseteq\IP^n$ left invariant by $X$. 
This is more directly achieved not by bounding the degree of a possible 
$V$, but 
by bounding the degrees of the polynomials necessary to cut out $V$. These 
degrees are bounded by the Castelnuovo-Mumford regularity of $V$, defined 
below.

Given a subscheme $V\subseteq\IP^n$, let $I$ be its saturated 
homogeneous ideal, and
\[
\cdots\to F_i\to\cdots\to F_0\to I\to 0
\]
the minimal graded free resolution of $I$. For each non-negative 
integer $i$, let $b_i$ be the maximum of the degrees of the generators of 
$F_i$. The 
Castelnuovo-Mumford regularity, or simply regularity, of $V$ is the integer 
$r$ given as
\[
r=\max(b_i-i\,|\,i=0,1,\dots)
\]
The above definition of regularity is of a rather arithmetical nature. 
However, if 
$V$ is arithmetically Cohen-Macaulay, by which we 
mean that its ring of homogeneous 
coordinates is Cohen-Macaulay, the regularity acquires a more geometric 
meaning. More precisely, cut 
out $V$ by as many general hyperplanes as its dimension to obtain a 
set $\Gamma$ of points. Then the regularity of $V$ is equal to that 
of $\Gamma$, and the regularity of $\Gamma$ is the smallest integer $r$ 
such that for each $P\in\Gamma$ there is a hypersurface of degree $r-1$ 
passing through all the points of $\Gamma$ but $P$. So, points disposed in 
special position tend to have higher regularity.

Complete intersections are arithmetically Cohen-Macaulay. The regularity 
of a hypersurface is its degree. Hence, the next theorem is a 
generalization of the inequality found by Cerveau and Lins Neto on the 
complex plane.

\begin{theorem} Let $C\subseteq\IP^n$ be a reduced, arithmetically 
Cohen-Macaulay curve with at most ordinary nodes for 
singularities, degree $d$ and Castelnuovo-Mumford regularity $r$. 
If $C$ is 
left invariant by a vector field $X$ of degree $m$ on $\IP^n$, 
and contains only finitely many singularities of $X$, then 
$r\leq m+2$, with equality only if $C$ is reducible or $p|d$.
\end{theorem}

A generalization of Theorem~1 is stated as Theorem 18, and 
follows from the main result of this note, Theorem 17. The assumption 
that $C\subseteq\IP^n$ be arithmetically Cohen-Macaulay is not so 
restrictive as it may seem. If $C$ is smooth, it's the same as assuming 
that $C$ be projectively normal. If $n=3$, it follows from 
the Hilbert-Burch Theorem that the set of arithmetically Cohen-Macaulay 
curves in $\IP^3$ is open in the Hilbert scheme; see \cite{Se}.

\begin{corollary} 
Let $C\subseteq\IP^n$ be a reduced curve 
with at most ordinary nodes for singularities. Assume $C$ is the 
complete intersection 
of hypersurfaces of degrees $d_1,d_2,\dots,d_{n-1}$. 
If $C$ is left invariant by a vector field $X$ of degree $m$ on $\IP^n$, 
and contains only finitely many singularities of $X$, then
\[
d_1+\cdots+d_{n-1}\leq m+n,
\]
with equality only if $C$ is reducible or $p|(d_1\cdots d_{n-1})$. 
\end{corollary}

As complete intersections are arithmetically Cohen-Macaulay, to 
get the above inequality it is enough to observe that the regularity of 
$C$ is $d_1+\cdots+d_{n-1}-n+2$, what follows from considering the 
Koszul resolution of the homogeneous ideal of $C$. For 
$C$ irreducible, our inequality is sharper than that given in \cite{CCF}. 
(The inequality given in \cite{CCF} is recovered more simply 
by our Proposition 12.)

Our technique is quite simple, and applicable in greater generality. 
Given a reduced, closed subscheme $V\subseteq\IP^n$ of dimension $s$, 
we can project it to a hypersurface $\ol V\subseteq\IP^{s+1}$. If the center 
of projection is general enough, then $\ol V$ will have 
the same degree of $V$ and singularities that are not 
much worse. (If $V$ is a curve with at most ordinary nodes for 
singularities, so is $\ol V$.) Now, if $X$ is a vector field on $\IP^n$ 
leaving $V$ invariant, can we find a vector field $\ol X$ on $\IP^{s+1}$ of 
degree not much higher than that of $X$ leaving $\ol V$ invariant? If so, 
we may use the bounds known for hypersurfaces. Well, the answer to our 
question is yes, if $V$ is arithmetically Cohen-Macaulay. As our Theorem~17 
shows, the difference between the degrees of $\ol X$ and $X$ need not 
be higher than $d-r$, where $d$ is the degree of $V$ and $r$ its 
regularity. To prove Theorem~17 we use elimination in a quite elementary 
way. 

Theorem~17 reduces our problem of bounding regularity to that of bounding 
the degree of a hypersurface. For curves, we can use the result by 
Cerveau and Lins Neto. This application of Theorem 17 is our Theorem 18, 
from which Theorem 1 above follows. 

For curves which are not arithmetically Cohen-Macaulay, 
the next theorem can be applied:

\begin{theorem} Let $C\subseteq\IP^n$ be a reduced curve with degree $d$ 
and at most ordinary nodes for singularities. 
Assume $C\subseteq Z\cap H$, where $Z$ is an arithmetically 
Cohen-Macaulay surface of degree $e$, Castelnuovo-Mumford regularity $r$, 
and $H$ is a hypersurface of degree~$f$ containing no irreducible 
component and just finitely many singularities of $Z$. If 
$C$ is left invariant by a vector field $X$ of degree $m$ on 
$\IP^n$, and contains only finitely many singularities of $X$, then
\[
d\leq m+f(e-1)-r+3,
\]
with equality only if $C$ is reducible or $p|d$.
\end{theorem}

The above theorem follows as well from our Theorem 18. 
If $C\subseteq\IP^n$ is a curve with at most ordinary nodes for 
singularities, it follows from the Bertini theorems of \cite{AK} that 
$C$ is contained in a 
smooth, complete intersection surface $Z\subseteq\IP^n$. So, the 
above theorem is always applicable. 

Using the result by Brunella and Mendes, we can generalize Theorem 1.

{\vskip0.2cm\parindent=0pt {\bf Theorem 1$^*$.} \emph{Let 
$V\subseteq\text{\bf P}^n_{\hskip-0.1cm\text{\bf C}}$ 
be a reduced, arithmetically Cohen-Macaulay subscheme of dimension $s$ and 
Castelnuovo-Mumford regularity $r$. Assume $V$ projects from a general 
center to a hypersurface with at most normal-crossings singularities. 
Let $X$ be a vector field of degree $m$ on 
$\text{\bf P}^n_{\hskip-0.1cm\text{\bf C}}$ leaving $V$ invariant. 
Assume each irreducible component of $V$ contains a non-singular point of 
$X$. Then $r\leq m+s+1$.}\vskip0.2cm}

Even when the subscheme 
$V\subseteq\text{\bf P}^n_{\hskip-0.1cm\text{\bf C}}$ is smooth, 
even when the center of projection is general, 
$V$ may project to a hypersurface with worse singularities than 
normal-crossings.

This note relies heavily on a few 
results of Commutative Algebra. For the benefit 
of the non-specialist, I tried to give precise references to these results 
as much as possible. 
I'm grateful to S. Kleiman, L.G. Mendes, J.V. Pereira, A. Simis, M. Soares 
and B. Ulrich for many helpful discussions. 
I thank also A. Garcia for working out with me Remark~14.

\section{Vector fields}
\label{section2}

\begin{subsct} \emph{Vector fields.} 
Let $V$ be a scheme over $k$ and $\mathcal L$ an 
invertible sheaf on $V$. A \emph{vector field on 
$V$ with coefficients in $\mathcal L$} is a map 
of $\Ocal_V$-modules $X\:\Omega^1_V\to\mathcal L$. 

Let $m$ be a non-negative integer. Using the Euler sequence,
\begin{equation}
\label{Euler}
0\to\Omega^1_{\IP^n}\to\Ocal_{\IP^n}(-1)^{n+1}\to\Ocal_{\IP^n}\to0,
\end{equation}
we see that a vector field $X\:\Omega^1_{\IP^n}\to\Ocal_{\IP^n}(m-1)$ 
is induced by a (homogeneous) vector field,
\[
G_0\d_0+\cdots+G_n\d_n,
\]
on $\text{\bf A}^{n+1}_k$, 
where $G_0,\dots,G_n\in k[t_0,\dots,t_n]$ are homogeneous of degree $m$, and 
$\d_0,\dots,\d_n$ are the standard derivations with respect 
to $t_0,\dots,t_n$. We say that \emph{$X$ is induced by $\sum G_i\d_i$}. 
Moreover, any two homogeneous vector fields on $\text{\bf A}^{n+1}_k$ 
determine the same vector field on $\IP^n$ if and only if they differ by a 
multiple of the \emph{radial} (or \emph{Euler}) vector field, which is 
induced by $\sum t_i\d_i$. By convention, $X$ has \emph{degree} $m$. 

The vector field $X$ induced by $\sum G_i\d_i$ may be 
viewed as the rational map sending a point $(t_0:\dots:t_n)\in\IP^n$ 
to $(G_0(t_0,\dots,t_n):\dots:G_n(t_0,\dots,t_n))$. The line 
passing by these two points is called the \emph{direction given by $X$ at 
$(t_0:\dots:t_n)$}. The rational map is defined away from the 
closed subscheme $V\subseteq\IP^n$ cut out by the maximal minors 
of the matrix 
\[
\left[\begin{matrix}
t_0&\hdots&t_n\\
G_0&\hdots&G_n
\end{matrix}\right].
\]
A point $P\in V$ is called \emph{a singularity of $X$}, or 
\emph{singular for $X$}.
\end{subsct}

\begin{subsct} \emph{Integral subschemes.} 
Given a scheme $V$ over $k$ and a vector field 
$X\:\Omega^1_V\to\mathcal L$, a subscheme $W\subseteq V$ is said to 
be \emph{left invariant by $X$}, or \emph{an integral subscheme of $X$}, 
if $X$ restricts to a vector field on $W$. 
More precisely, $W$ is left invariant by 
$X$ if there is a vector field $Y\:\Omega^1_W\to\mathcal L|_W$ 
making the following diagram commute:
\begin{equation}
\label{diagram}
\begin{CD}
\Omega^1_V|_W @>X|_W>> \mathcal L|_W\\
@VVV @|\\
\Omega^1_W @>Y>> \mathcal L|_W
\end{CD}
\end{equation}
where the left vertical map is the natural restriction. If 
\eqref{diagram} commutes, we say that $X$ \emph{restricts to} $Y$ or 
$Y$ \emph{lifts to} $X$.

Algebraically, if $X$ is the vector field on $\IP^n$ induced 
by $\sum G_i\d_i$, and $V\subseteq\IP^n$ is the closed subscheme cut out 
by the saturated, homogeneous ideal $I\subseteq k[t_0,\dots,t_n]$, 
then $X$ leaves $V$ 
invariant if and only if $\sum G_i\d_i F\in I$ for each $F\in I$. In fact, 
if $I$ is any homogeneous ideal cutting out $V$, then 
$\sum G_i\d_i F\in I$ for each $F\in I$ only if $V$ is left invariant by 
$X$. The converse does not hold in general, 
if $I$ is not saturated. For instance, the variety cut out by 
$t_0$ is left invariant by the 
vector field induced by $\d_1$, though 
$\d_1(t_0t_1)$ is not in the ideal generated by 
$t_0^2,t_0t_1,\dots,t_0t_n$.

Let $X$ be a vector field on $\IP^n$ of degree $m$. If $X\neq 0$, 
we can give a more geometric meaning 
to the degree. Indeed, let $H\subseteq\IP^n$ be a hyperplane 
not left invariant by $X$. The points of $H$ at which the direction given 
by $X$ lies in $H$ are those 
in the degeneracy scheme $Z$ of the natural map 
$\Omega^1_{\IP^n}|_H\to\Omega^1_H\oplus\Ocal_H(m-1)$ induced by $X$. 
Now, from the 
Euler sequences of $\IP^n$ and $H$ we get that the determinants of 
$\Omega^1_{\IP^n}$ and $\Omega^1_H$ are $\Ocal_{\IP^n}(-n-1)$ and 
$\Ocal_H(-n)$. So, $Z$ is the zero scheme of a section of $\Ocal_H(m)$, 
whence has degree $m$.
\end{subsct}

\section{Hypersurfaces and curves}
\label{section3}

\begin{subsct} \emph{Hypersurfaces.} 
Let $V\subseteq\IP^n$ be the hypersurface cut out by a 
homogeneous $F\in k[t_0,\dots,t_n]$. Then the vector field $X$ on 
$\IP^n$ induced by 
$\sum G_i\d_i$ leaves $V$ invariant if and only if
\begin{equation}
\label{invariant1}
G_0\d_0 F+\cdots+G_n\d_n F=PF
\end{equation}
for a certain $P\in k[t_0,\dots,t_n]$. We can easily construct vector 
fields that leave $V$ invariant, namely those induced by 
\begin{equation}
\label{all}
\sum_{i<j}P_{i,j}(\d_j F\d_i-\d_i F\d_j),
\end{equation}
with $P_{i,j}\in k[t_0,\dots,t_n]$ homogeneous of equal degrees. 
Conversely, the theorem below holds.
\end{subsct}

\begin{theorem} 
\label{thm7}
Let $V\subseteq\IP^n$ be the hypersurface cut out by $F\in k[t_0,\dots,t_n]$ 
homogeneous of degree $d$. 
If $V$ is smooth and 
$p{\not|}d$, then each vector field on $\IP^n$ 
leaving $V$ invariant is induced by \eqref{all} for certain homogeneous 
$P_{i,j}\in k[t_0,\dots,t_n]$ of equal degrees.
\end{theorem}

\begin{proof} The proof we present below is well-known. 
The argument appeared in \cite{L},~Part~c~of~Example~7,~p.~892, where it was 
attributed to Zariski.

Let $X$ be a vector field on $\IP^n$ that leaves $V$ 
invariant, say induced by $\sum G_i\d_i$. Then \eqref{invariant1} 
holds for a certain $P\in k[t_0,\dots,t_n]$. Since $p{\not|}d$, 
Euler's formula,
\[
dF=t_0\d_0 F+\cdots+t_n\d_n F,
\]
implies that \eqref{invariant1} is equivalent to
\begin{equation}
\label{invariant2}
G'_0\d_0 F+\cdots+G'_n\d_n F=0,
\end{equation}
where $G'_i:=G_i-Pt_i/d$ for $i=0,\dots,n$. Note that 
$X$ is also induced by $\sum G'_i\d_i$.

Let $S:=k[t_0,\dots,t_n]$. Let $I\subseteq S$ 
be the ideal generated by the partial derivatives of 
$F$:
\[
I:=(\d_0 F,\dots,\d_n F)\subseteq S.
\]
As $p{\not|}d$, by Euler's formula, $F\in I$. 
Since $V$ is smooth, $\sqrt I=(t_0,\dots,t_n)$. Hence, 
the partial derivatives $\d_0 F,\dots,\d_n F$ form a system of parameters 
for the local ring $\wt S$, where $\wt S$ is the localization of $S$ 
at $(t_0,\dots,t_n)$. Now, since $\wt S$ is regular, it is Cohen-Macaulay 
by \cite{Ma}, Thm.~17.8, p.~137. So, the sequence $\d_0 F,\dots,\d_n F$ 
is regular on $\wt S$ by \cite{Ma}, Thm.~17.4, p.~135. It follows 
that $\d_0 F,\dots,\d_n F$ is regular on $S$ as well, whence 
its Koszul complex is an exact sequence by \cite{Ma}, Thm.~16.5, p.~128. 
In particular, the module of relations,
\begin{equation}
\label{relations}
R:=\{(P_0,\dots,P_n)\in S^{n+1}|\  P_0\d_0 F+\cdots+P_n\d_n F=0\},
\end{equation}
is generated by the trivial ones,
\[
\d_i Fe_j-\d_j Fe_i\quad\text{for $0\leq i,j\leq n$}.
\]
As $(G'_0,\dots,G'_n)\in R$ by \eqref{invariant2}, and 
$X$ is induced by $\sum G'_i\d_i$, the proof is finished.
\end{proof}

\begin{theorem}
Let $X$ be a non-zero vector field on $\IP^n$ of degree $m$. 
Let $V\subseteq\IP^n$ be a hypersurface of degree 
$d$. If $V$ is smooth, left invariant by $X$, and $p{\not|}d$, then 
$d\leq m+1$.
\end{theorem}

\begin{proof} By Theorem~7, 
the vector field
$X$ is induced by \eqref{all} for 
certain homogeneous $P_{i,j}\in k[t_0,\dots,t_n]$ of the appropriate 
degrees. 
Since $X\neq 0$, we have $P_{i,j}\neq 0$ for certain $i,j$. 
Since $X$ has degree $m$, we have $\deg P_{i,j}=m-d+1$. 
As $\deg P_{i,j}\geq 0$, we get $d\leq m+1$.
\end{proof}

\begin{remark} 
\label{rmk9}
In the proof of Theorem 7 we observed the 
following fact: If $V\subseteq\IP^n$ is the hypersurface cut out 
by a homogeneous $F\in k[t_0,\dots,t_n]$ 
of degree $d$, and $p{\not|}d$, then 
any vector field on $\IP^n$ leaving $V$ invariant 
is induced by $\sum G_i\d_i$ for $(G_0,\dots,G_n)\in R$, where $R$ is 
the module of relations of the partial derivatives of $F$, 
defined in \eqref{relations}. So, if $q$ is the minimum degree of a 
non-zero relation in $R$, then the minimum degree of a non-zero vector 
field on $\IP^n$ leaving $V$ invariant is $q$ as well. 
The number $q$ measures how singular $V$ is. 
We have $q\leq d-1$ with equality if $V$ is smooth, as seen in 
the proof of Theorem 7. On the other hand, the worst case $q=0$ occurs 
if (and only if for $p=0$) the hypersurface $V$ is a cone. In $\IP^2$ a 
cone is integral only if it is a line. The plane curve cut out by 
$F:=t_0^at_1^b-t_2^{a+b}$, for $a$ and $b$ coprime, is an example of an 
integral plane curve with $q=1$. These well-known examples show that the 
answer to Poincar\'e's question is negative in general.
\end{remark}

\begin{subsct} \emph{The dualizing sheaf.} 
If $C\subseteq\IP^n$ is a reduced curve with at most ordinary nodes 
for singularities, then $C$ is a local complete intersection. (Indeed, 
applying repeatedly Theorem~7 on p.~787 of \cite{AK}, we see that 
$C$ is a divisor of a smooth, complete intersection surface in $\IP^n$.) 
So, the dualizing sheaf of $C$ is invertible by 
\cite{H}, Thm.~7.11, p.~245.
\end{subsct}

\begin{proposition} Let $C\subseteq\IP^n$ be a reduced curve 
with at most ordinary nodes for singularities. Let 
$\w$ be the dualizing sheaf of $C$. Let $\mathcal L$ be an 
invertible sheaf on $C$ and $X\:\Omega^1_C\to\mathcal L$ a vector field.  
If $\w\otimes\mathcal L^{-1}$ is ample, then $X=0$.
\end{proposition}

\begin{proof} Let $\pi\:\wt C\to C$ be 
the normalization of $C$ and $\eta\:\Omega^1_C\to\pi_*\Omega^1_{\wt C}$ 
the adjoint to the induced pullback map. Since 
$C$ has at most ordinary nodes for singularities, we can check locally 
that $\eta$ is surjective. As $\pi$ is 
birational, the kernel of $\eta$ is the torsion sheaf of $\Omega^1_C$. 
So
\[
\text{Hom}(\Omega^1_C,\mathcal L)=
\text{Hom}(\pi_*\Omega^1_{\wt C},\mathcal L).
\]
By Serre's duality (\cite{H}, Ch. III, \S 7) on $C$, 
the projection formula, and Serre's duality on $\wt C$,
\begin{align*}
\text{Hom}(\pi_*\Omega^1_{\wt C},\mathcal L)=&
H^1(\pi_*\Omega^1_{\wt C}\otimes\w\otimes\mathcal L^{-1})^*\\
=&H^1(\Omega^1_{\wt C}\otimes\pi^*(\w\otimes\mathcal L^{-1}))^*\\
=&H^0(\pi^*(\mathcal L\otimes\w^{-1})).
\end{align*}
Now, as $\w\otimes\mathcal L^{-1}$ is ample on $C$, the degree 
of $\pi^*(\mathcal L\otimes\w^{-1})$ on each component of $\wt C$ 
is negative. Thus $H^0(\pi^*(\mathcal L\otimes\w^{-1}))=0$, 
and so $X=0$.
\end{proof}

\begin{proposition} 
\label{prop12}
Let $C\subseteq\IP^n$ be a reduced curve 
with at most ordinary nodes for singularities. Assume $C$ is the 
complete intersection 
of hypersurfaces of degrees $d_1,d_2,\dots,d_{n-1}$. 
Let $m$ be an integer and $X\:\Omega^1_C\to\Ocal_C(m-1)$ a vector 
field. If $m<d_1+\cdots+d_{n-1}-n$ then $X=0$.
\end{proposition}

\begin{proof} The normal bundle of $C$ in $\IP^n$ is 
$\Ocal_C(d_1)\oplus\cdots\oplus\Ocal_C(d_{n-1})$. By 
\cite{H}, Thm.~7.11, p.~245, the dualizing sheaf $\w$ of $C$ 
satisfies 
\[
\w\cong\Ocal_C(d_1+\cdots+d_{n-1}-n-1).
\]
Since $m<d_1+\cdots+d_{n-1}-n$, the sheaf $\w\otimes\Ocal_C(1-m)$ is 
ample, and hence $X=0$ by Proposition~11.
\end{proof}

\begin{theorem}
\label{thm13}
Let $X$ be a non-zero vector field on $\IP^2$ of degree $m$. 
Let $C\subseteq\IP^2$ be an integral curve of degree 
$d$ whose only singularities are ordinary nodes. If $C$ is 
left invariant by $X$ and $p{\not|}d$, then $d\leq m+1$.
\end{theorem}

\begin{proof} The theorem was proved by Cerveau and Lins Neto for 
$k=\text{\bf C}$; see \cite{CL}, Thm. 1 and the remark thereafter on p. 891. 
We will prove the 
more general statement above in a forthcoming work. In the next remark, 
we will show what happens when $p|d$.  
\end{proof}

\begin{remark} Let $C\subseteq\IP^2$ be the curve cut out by a 
homogeneous $F\in k[t_0,\dots,t_n]$ of degree~$d$ divisible by $p$. 
Then $C$ is left invariant by a non-zero vector field 
$X$ on $\IP^2$ of degree $d-2$. Moreover, $X$ is unique 
because $\text{Hom}(\Omega^1_{\IP^2},\Ocal_{\IP^2}(-3))=0$. In concrete 
terms, write 
\[
F=t_0H_0+t_1H_1+t_2H_2.
\]
Then $X$ is induced by $\sum G_i\d_i$, where $G_0,G_1,G_2$ are given by:
\[
G_0:=\d_1 H_2-\d_2 H_1,\quad 
G_1:=\d_2 H_0-\d_0 H_2,\quad 
G_2:=\d_0 H_1-\d_1 H_0.
\]

For a numerical example, 
if $F=t_2^{d-1}t_0+t_0^{d-1}t_1+t_1^{d-1}t_2$ then $X$ is induced by 
\begin{equation}
\label{jouanolou}
t_1^{d-2}\d_0+t_2^{d-2}\d_1+t_0^{d-2}\d_2.
\end{equation}
If $p=0$, the vector field induced 
by \eqref{jouanolou} is known as Jouanolou's vector field, 
and its most striking property is that it leaves no curve invariant for 
$d\geq 4$; see \cite{J}.
\end{remark}

\section{Regularity and central projections}
\label{section4}

\begin{subsct} \emph{Regularity.} 
Let $\mathcal F$ be a coherent sheaf on $\IP^n$. Given an integer 
$m$, we say that $\mathcal F$ is \emph{$m$-regular} if 
$H^i(\mathcal F(m-i))=0$ for each positive integer $i$. 
By \cite{H},~Thm.~5.2,~p.~228, there is an integer $m$ such 
that $\mathcal F$ is $m$-regular. If $r$ is minimum 
among the integers $m$ such that $\mathcal F$ is $m$-regular, we call 
$r$ the \emph{Castelnuovo-Mumford regularity}, 
or simply \emph{regularity}, of $\mathcal F$. 

Suppose $\mathcal F$ is $m$-regular. By the 
proposition on page 99 of \cite{Mu1}, the sheaf 
$\mathcal F(m)$ is generated by its 
global sections. So, if $\mathcal F$ is a non-zero sheaf of ideals, then 
$m\geq 0$, with equality only if $\mathcal F=\mathcal O_{\IP^n}$. 
Hence, the regularity of a non-zero coherent sheaf of ideals is well-defined 
and non-negative. By definition, the regularity of $\IP^n$ is 1, and 
the regularity of a proper, closed subscheme of $\IP^n$ is that of 
its sheaf of ideals. 

There is also an algebraic notion of regularity, given as follows; 
see \cite{E},~p.~505. Let $S:=k[t_0,\dots,t_n]$ and $M$ a non-zero, 
finitely generated graded $S$-module. By Hilbert Syzygy Theorem 
(\cite{E}, Cor. 19.7, p. 474), there is a minimal graded free resolution of 
$M$,
\[
0\to F_m\to\cdots\to F_i\to\cdots\to F_0\to M\to 0,
\]
of length $m\leq n+1$. For each $i=0,\dots,m$, let $b_i$ be the maximum of 
the degrees of the generators of 
$F_i$. Let $r:=\max(b_i-i\,|\,i=0,\dots,m)$. We say that $r$ is 
the \emph{regularity} of $M$. 
It's clear from the definition that if $J\subset S$ is a non-zero, proper 
homogeneous ideal with regularity $r$, then the regularity of $S/J$ is 
$r-1$. By convention, the regularity of the zero module is 1.

For the correspondence between the algebraic and geometric notions of 
regularity, see \cite{E},~Ex.~20.20,~p.~516. In particular, if 
$V\subseteq\IP^n$ is the closed subscheme cut out by a saturated, 
homogeneous ideal $J\subseteq S$, then the regularity of $V$ is equal to 
that of $J$. It follows that the regularity of $V$ is an upper 
bound for the degrees of the hypersurfaces necessary to cut out $V$. 

It is a conjecture by D. Eisenbud and S. Goto 
that the regularity of a non-degenerate, integral closed subscheme 
$V\subseteq\IP^n$ of degree $d$ and codimension $c$ is bounded above by 
$d-c+1$; see \cite{EG}, p. 93. 
When $V$ is a curve, this bound was proved by L. Gruson, R. Lazarsfeld and 
C. Peskine; see \cite{GLP}, Thm. 1.1, p. 494. 
In its full generality, the conjecture has not been proved yet.

Let $V\subseteq\IP^n$ be a closed subscheme of pure codimension $c$, 
and $I\subseteq S$ its saturated, homogeneous ideal. We say that $V$ is 
\emph{arithmetically Cohen-Macaulay} 
if $S/I$ is Cohen-Macaulay. If $V$ is non-degenerate, 
arithmetically Cohen-Macaulay, integral, and of degree $d$, 
its regularity is bounded above by $\lceil(d-1)/c\rceil+1$, 
which is sharper than the bound conjectured by Eisenbud and Goto; 
see \cite{T}, Cor. 2.3, p. 144.
\end{subsct}

\begin{subsct} 
\emph{Central projections.} Assume $n\geq 2$. Let $S:=k[t_0,\dots,t_n]$. 
Let $L\subseteq\IP^n$ be a linear 
subspace of dimension $\ell$, with $0\leq\ell\leq n-2$. Let 
$s_0,\dots,s_{n-\ell-1}\in S$ be linear forms cutting out $L$. 
We view $s_0,\dots,s_{n-\ell-1}$ as the homogeneous coordinates of 
$\IP^{n-\ell-1}$. Let $\ol S:=k[s_0,\dots,s_{n-\ell-1}]$. 

Let $V\subseteq\IP^n$ be the closed subscheme cut out by a saturated, 
homogeneous ideal $I\subseteq S$. Assume $L$ does not meet $V$. 
Then the projection with center $L$ gives rise to a finite 
map $V\to\IP^{n-\ell-1}$, whose image $\ol V\subseteq\IP^{n-\ell-1}$ 
is the closed subscheme cut out by the ideal $\ol I:=I\cap\ol S$. 
We call $\ol V$ the \emph{projection of $V$ with center $L$}.

Let $X$ be a vector field on $\IP^n$ induced by $\sum G_i\d_i$.
After a change of coordinates, we may 
assume $s_i=t_i$ for $i=0,\dots,n-\ell-1$. 
Suppose $X$ leaves $V$ invariant. Then
\begin{equation}
\label{proj}
G_0\d_0 F+\cdots+G_{n-\ell-1}\d_{n-\ell-1}F\in I
\end{equation}
for each $F\in\ol I$. 
Now, suppose there are 
$B\in S$ and $\ol G_0,\dots,\ol G_{n-\ell-1}\in\ol S$ homogeneous 
such that $BG_i\equiv\ol G_i\mod I$ for each $i=0,\dots,n-\ell-1$. 
Then, for each $F\in\ol I$, 
it follows from (\ref{proj}) that
\[
\ol G_0\d_0 F+\cdots+\ol G_{n-\ell-1}\d_{n-\ell-1}F\in\ol I.
\]
So the projection 
$\ol V$ is left invariant by the vector field $\ol X$ on $\IP^{n-\ell-1}$ 
induced by $\ol G_0\d_0+\cdots+\ol G_{n-\ell-1}\d_{n-\ell-1}$. If 
$X$ is of degree $m$, then $\ol X$ is of degree $m+\deg B$. If $V$ is 
reduced and $L$ is general, 
then the degree of $\ol V$ is equal to that of $V$, say $d$. If 
one could give a bound for $\deg B$ in terms of $d$ and $m$, then any 
bound on the degree of $\ol V$ in terms of the degree of $\ol X$ will 
give an inequality involving the degrees of $V$ and $X$.
\end{subsct}

\section{The main theorem}
\label{section5}

\begin{theorem} 
\label{Thm17}
Assume $n\geq 2$. 
Let $\ell$ be an integer such that $0\leq\ell\leq n-2$. Let 
$W\subseteq\IP^n$ be an arithmetically Cohen-Macaulay closed subscheme 
of degree $e$, Castelnuovo-Mumford regularity $r$ and pure 
codimension $\ell+2$. Assume 
\[
\wt W:=\{P\in W|\,\text{\rm codim}(T_{W,P},\IP^n)\leq\ell\}
\]
has codimension one in $W$. Let $V\subseteq W$ be a reduced, 
closed subscheme, and $X$ a vector field of degree $m$ on $\IP^n$ 
leaving $V$ invariant. 
Assume each irreducible component of $V$ contains a non-singular point 
of $X$. Then, for each general linear subspace $L\subseteq\IP^n$ of 
dimension $\ell$, there is a non-zero vector field of degree 
$m+e-r$ on $\IP^{n-\ell-1}$ leaving the projection of $V$ with 
center $L$ invariant. 
\end{theorem}

\begin{proof} Let $L$ and $L'$ be linear subspaces of $\IP^n$ of 
dimensions $\ell$ and $\ell+1$ satisfying $L\subseteq L'$. Assume 
$L$ and $L'$ are in general position. We claim that the following three 
conditions hold.
\begin{enumerate}
\item The intersection $W\cap L'$ is empty.
\item The projection of $W$ with center $L$ has degree $e$.
\item For each irreducible component $C\subseteq V$, the linear subspace 
of $\IP^n$ 
spanned by $L$ and a general $P\in C$ does not contain the direction given 
by $X$ at $P$.
\end{enumerate}
In fact, since 
$W$ has pure codimension $\ell+2$, a general linear subspace 
$L'\subseteq\IP^n$ of dimension $\ell+1$ does not meet $W$. Thus, 
Condition~1 holds. 

As for Condition~2, let $L''\subseteq\IP^n$ be a general linear 
subspace of dimension $\ell+2$ containing $L'$. 
As $W$ is of pure codimension $\ell+2$, 
the intersection $L''\cap W$ is finite, and contains at least one 
general point $P$ on each irreducible component $C\subseteq W$. 
Since $L$ is a general subspace of $L''$ of dimension $\ell$, we have 
$L\cap W=\emptyset$ and, for 
each irreducible component $C\subseteq W$, the linear subspace 
of $\IP^n$ 
spanned by $L$ and a general $P\in C$ intersects $W$ only at $P$. Now, 
as $\wt W$ has codimension 1 
in $W$, and $W$ is of pure dimension, for each irreducible component 
$C\subseteq W$, the general $P\in C$ satisfies 
$\text{\rm codim}(T_{W,P},\IP^n)>\ell$. Thus, since $L$ is general of 
dimension $\ell$, the linear subspace of $\IP^n$ spanned by $L$ and $P$ 
intersects $T_{W,P}$ only at $P$. So, for each irreducible 
component $C\subseteq W$ and each general $P\in C$, the linear subspace 
of $\IP^n$ spanned by $L$ and $P$ intersects $W$ only at $P$, and with 
multiplicity 1. Projecting $W$ with center $L$, we obtain a finite 
map which is an isomorphism over a dense open subscheme of the image. 
Hence, the image has degree equal to that of $W$, that is, degree $e$.

Finally, by hypothesis, for each irreducible 
component $C\subseteq V$, the general $P\in C$ is non-singular for $X$. 
Since $\ell\leq n-2$, and $L$ is general of dimension $\ell$, the 
direction given by $X$ at $P$ does not meet $L$. Our claim is 
proved.

Let $S:=k[t_0,\dots,t_n]$. 
After changing coordinates, we may assume $L'$ is cut out by 
$t_0,\dots,t_{n-\ell-2}$, and 
$L$ is cut out by $t_{n-\ell-1}$ in $L'$. 
As $W\cap L'=\emptyset$, 
there is a linear form in $k[t_0,\dots,t_{n-\ell-2}]$ 
which is not identically zero on any irreducible 
component of $W$. Up to subtracting from $t_{n-\ell-1}$ a certain multiple 
of this form, we may assume $t_{n-\ell-1}$ is not identically zero on 
any irreducible component of $W$. 
Let $\ol S:=k[t_0,\dots,t_{n-\ell-1}]$. 

Let $J\subseteq S$ 
denote the saturated, homogeneous ideal cutting out $W$, and 
put $J':=J+(t_0,\dots,t_{n-\ell-2})$. Let $\m:=(t_0,\dots,t_n)$. 
As $t_0,\dots,t_{n-\ell-2}$ cut out the empty set in $W$, we have 
$\sqrt{J'}=\m$. Then the quotient $S/J'$ is a vector space of 
finite dimension over $k$. By hypothesis, $S/J$ is 
Cohen-Macaulay of pure dimension $n-\ell-1$. So, $S_{\m}/J_{\m}$ is a 
local, Cohen-Macaulay ring of dimension $n-\ell-1$. 
As $t_0,\dots,t_{n-\ell-2}$ 
generate in $S_{\m}/J_{\m}$ an ideal of dimension zero, 
they form a regular sequence by \cite{Ma}, Thm. 17.4, p. 135. Now, the 
multiplicity of $S_{\m}/J_{\m}$ is the degree of $W$; see 
\cite{E}, Ex. 12.6, p. 276. Since $t_0,\dots,t_{n-\ell-2}$ form a 
regular sequence in $S_{\m}/J_{\m}$, 
by \cite{E}, Ex. 12.11, 
p. 279, the multiplicity of $S_{\m}/J_{\m}$ is equal to 
the dimension of 
$S_{\m}/J'_{\m}$ over $k$, which is the same 
as the dimension of $S/J'$. So $\dim_kS/J'=e$. 

Let $M_1,\dots,M_e$ be monomials in $S$ generating $S/J'$ over 
$k$. We may assume $M_1=1$. Let $w_i:=\deg M_i$ for $i=1,\dots,e$. By 
\cite{E},~Prop.~20.20,~p.~508, the regularity of $S/J'$ is equal to that 
of $S/J$, which is $r-1$. Since $S/J'$ has finite length, using 
\cite{E},~Ex.~20.18,~p.~516, we get
\begin{equation}
\label{maxw}
r=\max(w_i\,|\,i=1,\dots,e)+1.
\end{equation}
We may assume $w_e=r-1$. As $M_1,\dots,M_e$ generate $S/J'$, for each 
$i=1,\dots,e$ there are homogeneous $A_{i,j}\in k[t_0,\dots,t_{n-\ell-2}]$ 
for $j=1,\dots,e$ such that
\begin{equation}
\label{equiv}
t_{n-\ell-1}M_i\equiv A_{i,1}M_1+\cdots+A_{i,e}M_e \mod J.
\end{equation}
Let $A:=(A_{i,j})$ and $D:=\det(t_{n-\ell-1}-A)$. 
As $M_0=1$, and (\ref{equiv}) holds for each $i=1,\dots,e$, it follows 
that $D\in J$. Clearly, 
$D$ is monic in $t_{n-\ell-1}$, hence non-zero. In addition, $D$ is 
homogeneous of degree $e$, because the 
degree of $A_{i,j}$ is $1+w_i-w_j$ for each $i,j\in\{1,\dots,e\}$. 
As $D\in J\cap\ol S$, 
the projection $\ol W$ of $W$ with center $L$ is contained in the 
hypersurface cut out by $D$ in $\IP^{n-\ell-1}$. By Condition~2 above, also 
$\ol W$ has degree $e$. So, $D$ cuts out $\ol W$.

If $e=1$, also $r=1$. In this case, let $B:=1$. Then $B$ is non-zero, 
homogeneous of degree $e-r$, and the following property holds.
\begin{equation}
\label{mod}
\text{For each $G\in S$ there is $\ol G\in\ol S$ such that 
$BG\equiv\ol G\mod J$.}
\end{equation}

Suppose now that $e>1$. We claim there is $B\in\ol S$ non-zero, 
homogeneous of degree $e-r$, satisfying (\ref{mod}). Indeed, for 
each $i=1,\dots,e$, let $D_i$ be the determinant of the matrix obtained 
from $t_{n-\ell-1}-A$ by removing the $i$-th column and the last row. 
Put $B:=D_1$. For each $i=2,\dots,e$, let $A^{(i)}$ be 
the matrix obtained from $A$ by exchanging the last row 
with the row $R^{(i)}:=[R^{(i)}_1,\dots,R^{(i)}_e]$, where 
$R^{(i)}_j:=0$ for $j\in\{2,\dots,e\}-\{i\}$ and
\[
(R^{(i)}_1,R^{(i)}_i):=
\begin{cases}
(M_i,-1)&\quad\text{if $i<e$,}\\
(t_{n-\ell-1}M_e,0)&\quad\text{if $i=e$.}
\end{cases}
\]
By \eqref{equiv}, modulo $J$, the matrix $A^{(i)}$ has eigenvalue 
$t_{n-\ell-1}$, with $(M_1,\dots,M_{e-1},0)$ being one eigenvector 
if $i<e$ and $(M_1,\dots,M_{e-1},M_e)$ one if $i=e$. Since $M_1=1$, 
it follows that $\det(t_{n-\ell-1}-A^{(i)})\in J$. Expanding 
$\det(t_{n-\ell-1}-A^{(i)})$ by the cofactors of the last row, we 
get 
\begin{equation}
\label{det-equiv}
BM_i\pm D_i\pm t_{n-\ell-1}D_e\in J\quad\text{if $i<e$}
\end{equation}
and $t_{n-\ell-1}(BM_e\pm D_e)\in J$ if $i=e$. Now, as $t_{n-\ell-1}$ 
is not identically zero on any irreducible 
component of $W$, and $W$ is arithmetically Cohen-Macaulay, $t_{n-\ell-1}$ 
is regular in $S/J$. So $BM_e\pm D_e\in J$. As $D_1,\dots,D_e\in\ol S$, 
and $M_1,\dots,M_e$ generate $S/J'$, property (\ref{mod}) holds.

We show now that $B\not\in J$, whence $B\neq 0$. 
Indeed, if $B\in J$, as $BM_e\pm D_e\in J$, 
we have $D_e\in J$. Then, as (\ref{det-equiv}) holds, $D_i\in J$ for 
each $i=2,\dots,e-1$ as well. Since $D\neq 0$, there is $i\in\{1,\dots,e\}$ 
such that $D_i\neq 0$. Fix such $i$. As $J\cap\ol S$ is generated by $D$, 
and $D_i\in J$, the degree of $D_i$ would be at least $e$. However, as 
$\deg A_{e,i}>0$, the degree of $D_i$ is at most $e-1$. 
We have a contradiction. Hence $B\not\in J$. As $w_e=r-1$, we have 
$\deg A_{e,1}=r$, whence $\deg B=e-r$. So, $B$ is as claimed.

Let $\ol V$ be the projection of $V$ with center $L$. 
If $B$ cuts out a subscheme containing $\ol V$, 
then $\ol V$ is left invariant by the vector field on 
$\IP^{n-\ell-1}$ induced by $Bt_0^m\d_0$, which is non-zero because 
$B\neq 0$, and has degree $m+e-r$. The theorem is proved 
in this case.

Suppose now that $\ol V$ is not contained in the subscheme cut out by $B$. 
As $V$ is reduced, there is an irreducible component $C\subseteq V$ 
such that $B(P)\neq 0$ for the general $P\in C$. 
Let $G_0,\dots,G_n\in S$ homogeneous of degree $m$ such that 
$X$ is induced by $\sum G_i\d_i$. 
By (\ref{mod}), for each $i=0,\dots,n-\ell-1$, 
there is $\ol G_i\in\ol S$ such that $BG_i\equiv\ol G_i\mod J$. Let 
$\ol X$ be the vector field on $\IP^{n-\ell-1}$ induced by 
$\sum \ol G_i\d_i$. If 
$\ol X\neq 0$, then $\ol X$ has degree $m+e-r$. 
Moreover, since $X$ leaves $V$ invariant, $\ol X$ leaves 
$\ol V$ invariant; see 16. It remains 
to show that $\ol X\neq 0$. Suppose $\ol X=0$. Then there is 
$H\in\ol S$ such that $\ol G_i=Ht_i$, whence $BG_i\equiv Ht_i\mod J$, 
for each $i=0,\dots,n-\ell-1$. 
For the general $P\in C$ we have $B(P)\neq 0$; hence 
$G_i(P)t_j(P)=G_j(P)t_i(P)$ for each $i,j\in\{0,\dots,n-\ell-1\}$. So, 
the linear subspace of $\IP^n$ spanned by $L$ and $P$ contains the direction 
given by $X$ at $P$. We get a contradiction with Condition~3. 
So $\ol X\neq 0$. 
\end{proof}

\begin{proof} (Theorem 1$^*$) 
Let $e$ be the degree of $V$. Let 
$\ol V\subseteq\text{\bf P}^{s+1}_{\hskip-0.1cm\text{\bf C}}$ be 
the normal-crossings hypersurface obtained from $V$ by projection with 
general center. The degree of $\ol V$ is $e$ as well. 
By Theorem~17, the hypersurface $\ol V$ is left invariant by a 
non-zero vector field on $\text{\bf P}^{s+1}_{\hskip-0.1cm\text{\bf C}}$ 
of degree $m+e-r$. Then $e\leq m+e-r+s+1$ by \cite{BM}. So $r\leq m+s+1$.
\end{proof}

\begin{theorem} Let $C\subseteq\IP^n$ be a reduced curve with degree $d$ 
and at most ordinary nodes for singularities. 
Assume $C$ is contained 
in an arithmetically Cohen-Macaulay curve $C'$ of 
degree $e$ and Castelnuovo-Mumford regularity $r$ such that 
the tangent space $T_{C',P}$ of $C'$ at $P$ has dimension at most $2$ 
for all but finitely many points $P\in C'$. If 
$C$ is left invariant by a vector field $X$ of degree $m$ on 
$\IP^n$, and contains only finitely many singularities of $X$, then
\[
d\leq m+e-r+2,
\]
with equality only if $C$ is reducible or $p|d$. 
In particular, $r\leq m+2$ if $C=C'$.
\end{theorem}

\begin{proof} Suppose first that $n=2$. Then $r=e$. So, the desired 
inequality follows from Proposition~12 if $C$ is reducible or $p|d$, and 
Theorem~13 otherwise. 
We may thus assume $n\geq 3$. Let $\ol C\subseteq\IP^2$ be the 
projection of $C$ with center a general linear subspace $L\subseteq\IP^n$ 
of dimension $n-3$. Since $C$ has at most ordinary 
nodes for singularities, and $L$ is general, $\ol C$ has at most 
ordinary nodes as well. (The proofs given to \cite{H}, Prop.~3.5, p.~310 and 
Thm.~3.10, p.~313 apply.) By Theorem~17 with 
$V=C$ and $W=C'$, there is a non-zero vector field on $\IP^2$ of degree 
$m+e-r$ leaving $\ol C$ invariant. Now, apply Proposition~12 or Theorem~13 
to $\ol C$.
\end{proof}

\begin{proof} (Theorem 3) Let $C':=Z\cap H$. Since $Z$ is arithmetically 
Cohen-Macaulay, and $H$ does not contain any irreducible component of $Z$, 
also $C'$ is arithmetically Cohen-Macaulay by \cite{E}, Prop. 18.13, p. 455. 
It follows from \cite{E}, Ex. 12.11, p.~278 that the degree of 
$C'$ is $ef$. In addition, the regularity of $C'$ is $r+f-1$. (The 
proof given to \cite{E}, Prop. 20.20, p.~508 applies.) Since $Z$ is smooth 
at all but finitely many points of $C'$, the tangent space of $C'$ at all 
but finitely many points has dimension at most 2. 
From Theorem~18, we get
\[
d\leq m+ef-(r+f-1)+2,
\]
with equality only if $C$ is reducible or $p|d$.
\end{proof}

\begin{remark}
Let $V\subseteq\IP^n$ be a closed subscheme of pure codimension $e$. Then 
$V$ is contained in the complete intersection $W$ of 
$e$ hypersurfaces; see Cor. $4^*$ on p.~66 of \cite{Mu2}. The 
projection from a general center $L\subseteq\IP^n$ of dimension $e-2$ maps 
$V$ to a hypersurface of $\IP^{n-e+1}$. We may apply Theorem 17, 
as long as we verify that the required condition on the tangent spaces at 
points of $W$ holds. 

If $C\subseteq\IP^n$ is a reduced 
curve, it follows from a repeated application of 
Thm. 7 on p. 787 of \cite{AK} that $C$ is contained in a complete 
intersection surface $Z$ which is non-singular but at the singular points 
of $C$. Choosing a hypersurface $H\subseteq\IP^n$ containing $C$ but not 
$Z$, we may apply Theorem~3.
\end{remark}

\begin{example} Let $Q$ be a smooth quadric in $\IP^3$. Let 
$C\subseteq Q$ be an integral curve with bi-degree $(a,b)$ and at most 
ordinary nodes for singularities. Assume $p{\not |}(a+b)$. Assume $C$ is 
left invariant by a vector field $X$ on $\IP^3$ of degree $m$, and 
$C$ contains finitely many singularities of $X$. Applying Theorem~3, 
we get $m\geq\min(a,b)$. 
Equality may be achieved, and there is in general no bound for 
$\max(a,b)$ in terms of $m$; see the remark below. 
\end{example}

\begin{remark} Let $X$ be a vector field on $\IP^n$ 
of degree $m$. We might consider the following question, extending 
Poincar\'e's: Can we give a bound for the regularity of a curve in 
$\IP^n$ left invariant by $X$? As expected, 
the answer to this question is no. An example is given in \cite{CCF}: 
if $X$ is induced by $t_0t_2\d_0+t_1t_3\d_1$, for each positive 
integer $d$ the curve $C_d\subseteq\IP^n$ cut out 
by $t_2-dt_3,t_0t_2^{d-1}-t_1^d,t_4,\dots,t_n$ has regularity (and degree) 
$d$ and is left invariant by $X$. This example 
is quite ingeniously a variation of that given in Remark~9. Note that 
$X$ does not depend on $d$. 
Note however that $C_d$ is singular for $d\geq 3$. One 
might ask whether the answer to the question we raised 
above is yes if one considers only smooth curves. In the next 
paragraph, we give an example for $n=3$ to show that a bound cannot 
be found only in terms of $m$. This example may also 
be extended to one for each $n\geq 4$.

For each integer $d\geq 1$, let $\IP^1\to\IP^3$ be 
the map defined by taking $(t_0:t_1)$ to 
$(t_0^d:t_0^{d-1}t_1:t_0t_1^{d-1}:t_1^d)$, and 
$C_d$ its image. Then $C_d$ is non-singular, and 
is cut out in $\IP^3$ by 
\[
t_1t_2-t_0t_3,\quad t_1^{d-1}-t_2t_0^{d-2},\quad t_2^{d-1}-t_1t_3^{d-2}.
\]
Let $X_d$ be the vector field of degree 1 on $\IP^3$ induced by
\[
dt_0\d_0+(d-2)t_1\d_1-(d-2)t_2\d_2-dt_3\d_3.
\]
It's easy to check that $X_d$ leaves $C_d$ invariant for each positive 
integer $d$.
\end{remark}

\vspace{0.5cm}

{\smallsc Instituto de Matem\'atica Pura e Aplicada, 
Estrada Dona Castorina 110, 22460-320 Rio de Janeiro RJ, Brazil}

{\smallsl E-mail address: \small\verb?esteves@impa.br?}

\end{document}